\documentclass[12pt,reqno]{amsart}

\usepackage{amsmath,amssymb,amsfonts,amsthm,eucal,amscd,} 

\begin{document}

\numberwithin{equation}{section}

\newtheorem*{theorem*}{Theorem} 
\newtheorem{theorem}{Theorem} 
\newtheorem{lemma}{Lemma}[section] 
\newtheorem{proposition}[lemma]{Proposition} 
\newtheorem{corollary}[lemma]{Corollary} 
\theoremstyle{definition} 
\newtheorem{definition}[lemma]{Definition} 
\newtheorem{example}[lemma]{Example} 
\newtheorem{notation}[lemma]{Notation} 
\newtheorem{conjecture}[lemma]{Conjecture}
\newtheorem{remark}[lemma]{Remark} 
\theoremstyle{remark} 

\makeatletter
\@addtoreset{figure}{section}
\makeatother
\renewcommand{\thefigure}{\arabic{section}.\arabic{figure}}

\def\Rem#1{\noindent Remark~\ref{#1}}
\def\Def#1{\noindent Definition~\ref{#1}}
\def\Lem#1{\noindent Lemma~\ref{#1}} 
\def\Thm#1{\noindent Theorem~\ref{#1}} 
\def\Sec#1{\noindent Section~\ref{#1}} 
\def\Fig#1{\noindent Figure~\ref{#1}} 

\def\Prop{\noindent Proposition}
\def\Ex{\noindent Example}
\def\Cor{\noindent Corollary}

\long\def\BTHMm#1#2{\begin{theorem*}[#1]#2\end{theorem*}}
\long\def\BTHM#1#2{\begin{theorem}\LL{#1}#2\end{theorem}}
\long\def\BDF#1#2{\begin{definition}\LL{#1}#2\end{definition}}
\long\def\BPROP#1#2{\begin{proposition}\LL{#1}#2\end{proposition}}
\long\def\BLEM#1#2{\begin{lemma}\LL{#1}#2\end{lemma}}
\long\def\BNT#1#2{\begin{notation}\LL{#1}#2\end{notation}}
\long\def\BREM#1#2{\begin{remark}\LL{#1}#2\end{remark}}
\long\def\BEX#1#2#3{\begin{example}[#2]\LL{#1}#3\end{example}}
\long\def\BCOR#1{\begin{corollary}#1\end{corollary}}

\long\def\BEQ#1#2{\begin{equation}\LL{#1}#2\end{equation}}
\long\def\BB#1#2#3{\begin{#1}\LL{#2}#3\end{#1}} 
\long\def\BBF#1#2#3{\begin{figure}[h]#3\caption{#2}\LL{#1}\end{figure}} 

\def\Pr{\begin{proof}}
\def\rP{\end{proof}} 
\def\bydef#1{{\itshape\bfseries\rmfamily#1\,}}
\def\bydef#1{#1}

\def\i{iii}
\def\p{\mathcal P}
\def\X{X}
\def\A{\mathcal A}
\def\E{\mathcal E}
\def\J{\mathcal J}
\def\dd{\mathcal D}
\def\varpsi{\varPsi}
\def\uv{U}
\def\ud{D}
\def\udd{D^*}
\def\h{\mathcal H}
\def\L{\mathcal L}
\def\C{\mathbb C}
\def\R{\mathbb R}
\def\fone{\xi_R^{\phantom{y}}}
\def\ftwo#1{\psi_{R,#1}^{\phantom{y}}}
\def\zetaa#1#2{\zeta_{#2,#1}^{\phantom{y}}}

\def\bb#1#2{\begin{#1}\LL{#2}} 
\def\B#1{$\discretionary{}{}{}$#1$\discretionary{}{}{}$}
\newcounter{tmp}
\def\blist{\begin{list}{\hskip-19pt(\arabic{tmp})}{\usecounter{tmp}}}
\def\sitem{\parskip\par\item}

\def\boxx#1{\leavevmode\vbox{\hbox to 0pt{\hss\raise1.8ex\vbox 
to 0pt{\vss\hrule\hbox{\vrule\kern.75pt\vbox{\kern.75pt\hbox{\tiny #1}\kern.75pt}\kern.75pt\vrule}\hrule}}}\relax} 

\def\LL#1{\label{#1}\protect\boxx{#1}} 
\def\LL#1{\label{#1}}

\def\e#1{(\reF{#1})}
\def\reF#1{\ref{#1}}

\newenvironment{procLaim}{\it}{\par\smallskip}
\def\proclaim#1{\par\smallskip\noindent {#1}\bb{procLaim}{?}}
\def\endproclaim{\end{procLaim}} 

\def\zf{zeta function } 
\def\szf{spectral zeta function } 
\def\zfoap{zeta function of a polynomial} 
\def\zfotp{zeta function of the polynomial } 
\def\Str{Stric\-hartz} 
\def\Edef{\overset{\textup{\textsf{\tiny def}}}=}
\def\Si{Sier\-pi\'n\-ski} 
\def\Sig{Sier\-pi\'n\-ski gasket} 
\def\Sil{Sierpi\'nski lattice}
\def\iSig{in\-fi\-ni\-te Sierpi\'nski gasket} 
\def\Sip{Sierpi\'nski pre-gasket}
\def\Lp{Lap\-la\-ci\-an}
\def\bfit#1{{\bfseries\itshape#1\,}}

\def\<{\langle}
\def\>{\rangle}

\def\iiii{}
\def\iii{}
\def\ii{}
\def\ih{}
\def\ps{}
\def\pms{}
\def\lala{}
\def\ms{}
\def\sm#1{\mbox{$#1$}}\def\sm{}
\def\lg#1{\mbox{$#1$}}
\def\Lg#1{\mbox{$#1$}}
\def\LG#1{\mbox{$#1$}}

\def\split{}
\def\endsplit{}
\def\ndt{\noindent}\def\ndt{}
\def\np{\newpage}
\def\pb{\pagebreak}

\long\def\B{\par}

\def\<{\langle} 
\def\>{\rangle} 
\def\({\<\<} 
\def\){\>\>}
\def\tri{| \hskip-0.02in |\hskip-0.02in |}
\def\e{\alpha}
\def\p{\mathcal P}
\def\E{\mathcal E}
\def\J{\mathcal J}
\def\dd{\mathcal D}\def\D{\mathcal D}
\def\ddD#1{#1,#1}
\def\h{\mathcal H}
\def\Z{\mathcal Z}
\def\J{\mathcal J}
\def\A{\mathcal A}
\def\L{\mathcal Z}
\def\M{\mathcal M}
\def\N{\mathcal N}
\def\C{\mathbb C}
\def\R{\mathbb R}

\def\Edef{\overset{\textup{\textsf{ def}}}=}
\def\Si{{Sierpi\'nski}} 
\def\Sig{Sierpi\'nski gasket} 
\def\Sil{Sierpi\'nski lattice}
\def\iSig{infinite Sierpi\'nski gasket} 
\def\Sip{Sierpi\'nski pre-gasket}
\def\Lp{Laplacian}
\def\uv{U}
\def\ud{D}
\def\udd{D^ *}
\def\X{
X}
\def\bX{\bar 
X}

\def\dim{\textup{\textsf{dim}}}
\def\supp{\textup{\textsf{supp}}}
\def\mult{\textup{\textsf{mult}}}
\def\deg#1{\textup{\textsf{deg}}_\ii {\raise-.2ex\hbox{$_\ii #1$}}}
\def\degg{\textup{\textsf{deg}}}
\def\ind#1{_\ii {\raise-.1ex\hbox{$_\ii #1$}}}
\def\inds#1{
{\raise-.25ex
\hbox{$
_\ii 
#1$}}}

\def\z{^\ii {\text{\bfseries (0)}}}
\def\ellz{\ell^\ii {2\text{$(0)$}}}

\def\1{\mbox{$1$%
\hskip-.38em\rule{.44em}{.05ex}%
\hskip-.18em\rule{.06em}{1.5ex}%
\hskip-.16em{\raise1.48ex\hbox{\rule{.27em}{.05ex}}}\hskip0.1em%
}}

\long\def\NEWSECTION#1#2#3{\section{{#2}}
\B\centerline{#3}\B\B}
\long\def\NEWSECTION#1#2#3{\section{{#2}}\LL{#3}}

\def\SMASH{\smash}
\def\SMASH#1{\,{#1}\,}
\def\SMASH#1{\smash{#1}\vphantom{\frac12}}
\def\SMASH#1{\smash{#1}}
\def\SMASH{}

\title[Spectral zeta functions of fractals and polynomials]
{Spectral zeta functions of fractals and the complex dynamics of polynomials}
\author{Alexander Teplyaev}
\thanks{Research supported in part by  NSF grants DMS-0071575 and DMS-0505622}
\address{\noindent Department of Mathematics, 
University of Connecticut, Storrs CT 06269 USA}
\email{teplyaev@math.uconn.edu}
\date{\today}
\begin{abstract}\noindent 
We obtain formulas for the spectral zeta function 
of the \Lp\ on symmetric finitely ramified fractals, 
such as the Sier\-pi\'n\-ski gasket, and 
a fractal \Lp\ on the interval. 
These formulas contain 
a new type of zeta function associated 
with a polynomial (rational functions also can appear in this context). It is proved that 
this zeta function has a meromorphic 
continuation to a half plain 
with poles contained in an arithmetic progression. 
It is shown as an example that 
the Riemann zeta function is 
the zeta functions of a quadratic polynomial, which is 
associated with the Laplacian on an interval. 
The spectral zeta function of the \Sig\ is a product 
of the zeta function of a polynomial and a geometric part; 
 the poles of the former 
are canceled by the zeros of the latter. 
A similar product structure was 
discovered by M.L.~Lapidus for self-similar fractal strings. 
{\tableofcontents}
\end{abstract} 
\maketitle

\section{Introduction}

In this paper we obtain formulas for the 
spectral zeta function of the \Lp\ on 
symmetric finitely ramified (or p.c.f., see \cite{Ki1,Ki}) fractals, such as 
the \Sig\ and an interval with a fractal measure. 
These formulas involve a new type of zeta function associated 
with a polynomial, and we show that 
this zeta function has a meromorphic 
continuation to a strip beyond the half plane where it is analytic. 
In the case of the \Sig\ 
the spectral zeta function has a product structure, with one term 
which can be described as a ``geometric'' part, and the other term 
which is the \zf of a certain quadratic polynomial. Both of these 
features, the existence of a meromorphic continuation and a product 
formula, are reminiscent of the theory of one dimensional fractal strings 
first defined and studied by M.~Lapidus in \cite{L92}. We present basic definitions and 
results in \Sec{s2}, and the reader can find more details and references in \cite{LvF,LvF2,LM1,LM2,LP1,LP2}. 
In particular, equations (5.2)-(5.3) in \cite{L92} 
(see also section 1.3 of \cite{LvF}) contain 
the product formula 
\BEQ{e00}{\zeta_L(s)=
\zeta(s)\zeta_\mathcal L(s),} 
where 
$\SMASH{\zeta(s)}$ is 
the Riemann zeta function, 
 and 
$$ 
\SMASH{\zeta_\mathcal L(s)=\sum l_k^{ {s} }} 
$$ 
is the \bydef{geometric zeta function} of a 
fractal string $\mathcal L$, 
a disjoint collection 
of intervals of lengths $\SMASH{l_k}$. 
The \bydef{spectral zeta function} 
$\zeta_L(s)$ of a normalized Dirichlet Laplacian 
$L$ on $\mathcal L$ is defined as 
\BEQ{e-szf}{ 
\zeta_L(s)= 
\sum\lambda_j^{-{ {s} }\sm{/2}},}
where $\SMASH{{\lambda_j}}$ are the nonzero 
eigenvalues of $L$.

The first step in our work is to define and study a 
new type of zeta function 
associated with a polynomial (see \Sec{s3}). 
Suppose that $R(z)$ is a polynomial of 
degree $N$ such that 
$R(0)\pms=0$, 
$c\pms=R'(0)\pms>1$, and the Julia set of $R(z)$ 
lies in $\mathbb R_+$. 
Then we define the \bydef{zeta function of the polynomial} 
$R(z)$ 
as 
\BEQ{e0A}{
\SMASH{
\ 
\zetaa{z_0}   {R}(s)\ps=
\lim\limits_{n\to\infty}
\sum_{z\in R^{-n}\{z_0\}}
(c^n z)^{- {s} /2}
\ }} 
for 
\BEQ{e0B}{\ \text{Re}(s)\ps>d_R\ps=\frac{2\log N}{\log c}.} 
where $R^{-n}\{z_0\}$ means the set of $N^n$ preimages of $z_0$ under the $n$th composition power of the polynomial $R$. 
Moreover, we prove that there is an increasing unbounded sequence of positive real numbers $\lambda_j$ such that 
$\zetaa{z_0}   {R}$ can be represented by the Dirichlet series 
\BEQ{e0Alambda}{
\zetaa{z_0}   {R}(s)\ps=
\sum_{j=1}^\infty
\lambda_j^{- {s} /2}
} 
if (\ref{e0B}) holds. 
We give examples where $\lambda_j$ are actually eigenvalues of a \Lp, but it is possible to define $\lambda_j$
in terms of $R$ and $z_0$ only, without reference to any \Lp. 

Although the definition above make sense under less 
restrictive conditions on 
$R(z)$, 
they are automatically satisfied if the polynomial is 
associated with a self-similar \Lp\ (see \cite{Sh,MT,T01}), 
as it will be the case 
in the applications we consider. 

Our main result is that 
$\zetaa{z_0}   {R}(s)$ 
has a meromorphic continuation of the 
form 
\BEQ{e0}{ 
\zetaa{z_0}   {R}(s)=\frac{\fone(s)}{1-Nc^{- {s} /2}}+\ftwo{z_0}(s)
}
where $\fone(s)$ is analytic in $\mathbb C$ 
and $\ftwo{z_0}(s)$ is analytic for 
$ \text{Re}(s)\pms>0$. 
If $\mathcal J_R$ is totally disconnected, 
then there exists $\varepsilon\pms>0$ such that 
$\ftwo{z_0}(s)$ is analytic for 
$ \text{Re}(s)\pms> -\varepsilon$. 
The set of poles of $\zetaa{z_0}   {R}(s)$ is contained in 
$ 
\left\{\frac{2\log N+4in\pi}{\log c}\pms: n\in\mathbb Z\right\}
$, and there always is a pole at $d_R$. 

Independently of our work this result was recently 
(after the resent paper was completed) improved in \cite{DGV} based 
on the ideas of an earlier paper \cite{Grabner} by Peter Grabner. In particular, it is shown 
in  \cite{DGV} that $\zetaa{z_0} {R}$ actually has a meromorphic continuation to the 
entire complex plain for any polynomial $R$ satisfying the same conditions as in our paper. 
Moreover, \cite{DGV} contains a detailed analysis of 
$\zetaa{z_0}{R}$  in the case of the \Sig. 
However, an advantage of our method  is that it is 
 also applicable in the case when $R$ is a rational function satisfying 
appropriate assumptions. Therfore it is likely that our results can be extended to 
the class of fractals considered in \cite{MT}.

One of the reasons this new class of zeta functions is interesting 
is that the 
Riemann zeta function 
$\zeta(s)$ is the \zf of a quadratic polynomial (\Thm{thmRiemann}). 
This result is obtained by comparing two representations of the 
{Riemann} zeta function: as the Dirichlet series and as 
the \szf of the standard \Lp\ on an interval. 
Then the interval can be considered as a self-similar fractal, and 
the spectrum of the \Lp\ can be described in terms of the 
quadratic polynomial. This construction is not unique in the sense that 
 the 
Riemann zeta function can be represented as the \zfoap\ of any degree 
 (see \Thm{thmRiemann} for an example). 
 In \Sec{s4} this construction is generalized for 
fractal \Lp s on the interval, 
that is, \Lp s associated with a fractal 
self-similar measure.

In \Sec{s5} we obtain the following 
formula for the spectral zeta function of 
the \Lp\ $\Delta_\mu$ on the \Sig\ 
\begin{equation}\label{e1}
\zeta{\ind{{\Delta_\mu}}}(s)=
\tfrac12\zetaa{\frac34}   {R}(s)
\left(
\tfrac{1}{5^{ {s} /2}-3}+\tfrac{3}{5^{ {s} /2}-1}
\right)
+
\tfrac12\zetaa{\frac54}   {R}(s)
\left(
\tfrac{3\cdot5^{- {s} /2}}{5^{ {s} /2}-3}-\tfrac{5^{- {s} /2}}{5^{ {s} /2}-1}
\right)
\end{equation}
where $R(z)=z(5-4z)$. 

Combining (\ref{e0}) and (\ref{e1}) we obtain that 
the spectral 
zeta function of the \Sig\ has a form 
which resembles the spectral zeta function of a 
self-similar 
fractal string 
(\ref{e00}). For example, for the Cantor self-similar 
fractal string $\mathcal L$, that is 
the complement of the middle third Cantor set in $[0,1]$, we have 
$$
 \SMASH{\zeta_L(s)=\zeta(s)/(3^{\sm{\sm{\sm{}}} {s} }-2 )} .
$$ 
Thus the \szf of the \Sig\ has the same structure as 
the \szf of a disjoint collection 
of ``fractal intervals''.
However, this analogy is somewhat superficial because 
such ``fractal intervals'' so far have not been constructed or defined. 
This means that the spectrum has a product structure which does not manifest itself 
in the structure of the underling space. 

The poles of a \szf are called the \bydef{complex spectral dimensions} 
(see \cite{L92,LvF,LvF2}). They are $1$ and 
$\big\{ \frac{\log 2+2in\pi}{\log 3}: 
n\in\mathbb Z \big\}$ 
for the Cantor self-similar 
fractal string, and 
${\big\{\frac{2in\pi}{\log 5},\frac{\log 9+2in\pi}{\log 5}: 
n\in\mathbb Z\big\}}$ 
for the \Sig, see Figures~\ref{figCanCSD} and~\ref{figSigCSD}. 
An interesting feature of the product formula (\ref{e1}) is that, 
in the case of the \Sig, 
the poles of the 
``polynomial'' part of 
$\zeta{\ind{{\Delta_\mu}}}(s)$ 
are canceled by the 
zeros of the 
``geometric'' part of 
$\zeta{\ind{{\Delta_\mu}}}(s)$. 
We conjecture that such a cancelation takes place for 
other symmetric fractals, and 
 can be 
explained by the interplay between geometric and analytic properties of the fractals.

The main results of this paper were announced, without proofs, in \cite{T04} with a few typos, which we correct in the present paper.

\subsubsection*{Acknowledgment} 
The author is deeply grateful to M.~Lapidus for mentoring, 
encouragement and support, and to R.~Stric\-hartz 
and J.~Ki\-ga\-mi 
for many helpful discussions. 
The author very much appreciates interesting and useful remarks and references by P.~Grabner. The author thanks the anonymous
  referee for many helpful suggestions.

\NEWSECTION1{The spectral zeta function of a self-similar string} 
{s2} 

In this expository section we follow the work of Lapidus {\itshape et al.} 
(see \cite{L92,LvF,LvF2,LM1,LM2,LP1,LP2} and the references therein). 

\BDF{df01}{If $A$ is a nonnegative self-adjoint operator 
with discrete spectrum and positive eigenvalues $\{\lambda_j\}$, 
counted with 
multiplicities, 
then its \bydef{spectral zeta function} is 
$$
\zeta_A(s)=\sum\lambda_j^{-{ {s} }\sm{/2}}.
$$
Note that  zero  is not included in the sum, even if 
it is an eigenvalue of $A$.}

\BEX{ex02}{}{
If $A$={$-\frac{d^2}{dx^2}$} is a Neumann or Dirichlet Laplacian on 
$[0,1]$ then 
$$
\zeta_A(s)=
\sum_{j=1}^\infty\left(\frac{j^2}{\pi^2}\right)^{\hskip-7pt-{ {s} }\sm{/2}}
\hskip-7pt=
\ \pi^{ {s} }\zeta(s)
$$
where $\zeta(s)$ is the Riemann zeta function. The series converges for 
$\text{Re}(s)>1$, and has a meromorphic continuation to $\mathbb C$ with one 
pole at $s\ps=1$. 
}

\BDF{df03}{
The poles of $\zeta_A(s)$ are called the 
\bydef{complex spectral dimensions} of $A$.}

\BDF{df04}{
If $\mathcal L$ is a fractal string, that is, a disjoint collection 
of intervals of lengths $l_k$, 
then its \bydef{geometric zeta function} is 
$$
\zeta_\mathcal L(s)=\sum l_k^{ {s} }.
$$}

The next theorem is due to M.L. Lapidus \cite{L92,LvF,LvF2,LM1,LM2,LP1,LP2}. 
\BTHMm{}{If $A\ps=-\frac{d^2}{dx^2}$ is a 
Neumann or Dirichlet Laplacian on 
$\mathcal L$ then 
$$
\zeta_A(s)=\pi^{ {s} }\zeta(s)\zeta_\mathcal L(s).
$$}

\BEX{ex05}{a self-similar fractal string}{
Let $r\sm>1$ and $l_j=cr^{-k}$ for $N^k$ intervals, 
where $N$ is a positive 
integer, $k=0,1,2,...$ . Then 
$$
\zeta_\mathcal L(s)=\frac{c^{ {s} }}{1-Nr^{- {s} }},
\ \ \ \ 
\zeta_A(s)=\zeta(s)\frac{c^{ {s} }\pi^{ {s} }}{1-Nr^{- {s} }}
$$
If $N< r$ then the total length of $\mathcal L$ is finite 
and the spectral dimension is $1$, otherwise it is 
$\frac{\log N}{\log r}$. 
The 
complex spectral dimensions are $1$ and $\{\frac{\log N+2in\pi}{\log r}: n\in\mathbb Z\}$.

For example, for the Cantor self-similar 
fractal string, that is 
the complement of the middle third Cantor set in $[0,1]$, we have $c=\frac13$, $r=3$, $N=2$. 

\BBF{figCan}
{Cantor self-similar 
fractal string.}
{\begin{picture}(243,30)(0,-5)\setlength{\unitlength}{0.5pt}\small
\put(0,0){\line(1,0){486}}
\put(-7,9){$0$}
\put(157,17){$\frac13$}
\put(319,17){$\frac23$}
\put(486,9){$1$}
\linethickness{3.5pt}%
\setlength{\unitlength}{81pt}
\put(1,0){\line(1,0){1}}
\put(0,0){\setlength{\unitlength}{27pt}\put(1,0){\line(1,0){1}}
\put(0,0){\setlength{\unitlength}{9pt}\put(1,0){\line(1,0){1}}
\put(0,0){\setlength{\unitlength}{3pt}\put(1,0){\line(1,0){1}}}
\put(2,0){\setlength{\unitlength}{3pt}\put(1,0){\line(1,0){1}}}
}
\put(2,0){\setlength{\unitlength}{9pt}\put(1,0){\line(1,0){1}}
\put(0,0){\setlength{\unitlength}{3pt}\put(1,0){\line(1,0){1}}}
\put(2,0){\setlength{\unitlength}{3pt}\put(1,0){\line(1,0){1}}}
}
}
\put(2,0){\setlength{\unitlength}{27pt}\put(1,0){\line(1,0){1}}
\put(0,0){\setlength{\unitlength}{9pt}\put(1,0){\line(1,0){1}}
\put(0,0){\setlength{\unitlength}{3pt}\put(1,0){\line(1,0){1}}}
\put(2,0){\setlength{\unitlength}{3pt}\put(1,0){\line(1,0){1}}}
}
\put(2,0){\setlength{\unitlength}{9pt}\put(1,0){\line(1,0){1}}
\put(0,0){\setlength{\unitlength}{3pt}\put(1,0){\line(1,0){1}}}
\put(2,0){\setlength{\unitlength}{3pt}\put(1,0){\line(1,0){1}}}
}
}
\end{picture}}

\BBF{figCanCSD}
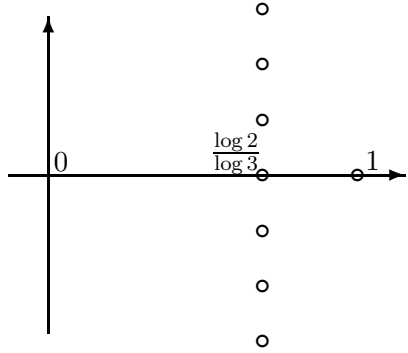
{Complex spectral dimensions of the 
Cantor self-similar 
fractal string.} 
{\begin{picture}(150,150)(-15,-70)\setlength{\unitlength}{0.5pt}\small 
\thicklines\setlength{\unitlength}{0.3pt}
\put(-50,0){\vector(1,0){500}}
\put(0,-200){\vector(0,1){400}}
\put(7,5){{$0$}}
\put(205,19){$\sm{\frac{\log 2}{\log 3}}$}
\put(400,7){${1}$}
\multiput(270,-210)(0,70){7}{\circle{12}}
\put(390,0){\circle{12}}
\end{picture}}
}

\NEWSECTION2{The zeta function of a polynomial}{s3}
\nopagebreak
We suppose that $R(z)$ is a polynomial of 
degree $N$ 
such that 
$R(0)\pms=0$ 
and 
$c\pms=R'(0)\ps>1$. We also assume that the Julia set $\mathcal J_R$ of 
$R(z)$ satisfies $\mathcal J_R\subset\mathbb R_+$. 
These assumptions imply immediately that the coefficients of $R$ are real. 
Since zero is a repulsive fixed point, its preimages are dense in 
 $J_R$ by the classical paper \cite{Br}. Thus, 
 $\mathcal J_R\subset\mathbb R_+$ 
if and only if there is a real interval $[0,a]\subset\mathbb R_+$  
such that $R^{-1}[0,a]\subset[0,a]$. 

By \cite{Br}, the Julia set in our 
situation is an interval or a Cantor set of Lebesgue measure zero. We 
denote the convex hull of $\mathcal J_R$ by $I$, which is an interval. 
By the same paper, if $\mathcal J_R$ is an interval, then $R^{-n}(I)= I$. 
Otherwise $R^{-n}(I)\varsubsetneq I$ consists of $N^n$ disjoint intervals. 
In both cases $R^{-n-1}(I)\subset R^{-n}(I)$ and 
$$J_R=\bigcap_{n=1}^\infty R^{-n}(I),$$ 
where $R^{-n}(\cdot)$ denotes the preimage under the $n$th 
composition (iteration) power of the polynomial $R(z)$. 
Also one has a unique decomposition $I=\bigcup_{j=1}^NI_j$, where 
$I_j$ are intervals with disjoint interior, and $R(z)$ is monotone 
on each $I_j$. 

In what follows we denote $I_0=\bigcap_{j=1}^NR(I_j)$. Then we have $R^{-1}(I_0)\subseteq I\subseteq I_0$.

\BDF{dfZFOAP}{The \bydef{zeta function of the polynomial} $R(z)$ is
\begin{equation}\label{ePZF}
\ 
\zetaa{z_0}   {R}(s)\ps=
\lim\limits_{n\to\infty}
\sum_{z\in R^{-n}\{z_0\}\backslash0}
(c^n z)^{- {s} /2}
\end{equation}
where $z_0\in I_0$, and $\text{Re}(s)\ps>d_R$. 
Here we define \BEQ{e-dR}{d_R\ps=\frac{2\log N}{\log c}.} 
}

\BTHM{propPZF}{For 
$\text{Re}(s)\ps>d_R\ps$ the limit (\ref{ePZF}) exists and is an analytic function of $s$.
Moreover, there is a  
nondecreasing unbounded sequence of positive real numbers $\lambda_j$ such that
\BEQ{e0Alambda-t}{
\SMASH{
\ 
\zetaa{z_0}   {R}(s)\ps=
\sum_{j=1}^\infty
\lambda_j^{- {s} /2}
\ }} 
where $\lambda_j$ are defined in terms of $R$ and $z_0$. 
}

\Pr 
By the monotonicity, the restriction of the polynomial $R(z)$ to each $I_j$ has an 
inverse, which we denote by $R^{-1}_j$. 
We assume that the intervals $I_j$ are numbered 
consecutively starting  from the origin. 
In particular, $R^{-n}_1(\cdot)$ is 
the $n$th composition power of 
the branch of the 
inverse function of $R(z)$ that passes through the origin.

We denote by 
$\mathcal J_+$ the part of $\mathcal J_R$ that is not in 
$R^{-1}_1(I)$, that is $$\mathcal J_+=\mathcal J_R \bigcap
\big(\bigcup_{k=2}^N R^{-1}_k(I)\big)= 
\bigcup_{k=2}^N R^{-1}_k(\mathcal J_R).$$ 

We denote by 
$\mathcal R$ the limit 
\BEQ{emathcalR}{\mathcal R(z) = \lim\limits_\ii {n\to\infty} c^\ii n R^\ii{-n}_\ii{1}(z).} 
Since $R'(0)=c$, this limit exists for any $z\in I_0$ and satisfies 
$\mathcal R(z) = c \mathcal R(R^{-1}_1(z))$ and is an increasing continuous function. 
The function $\mathcal R$ was first defined and studied by P.~Fatou in \cite{F}. Its analytic properties are well studied in complex dynamics, but we need only the elementary facts 
that 
the limit (\ref{emathcalR}) 
 exists and is an increasing continuous function for real $z\in I_0$. To show this, notice that $R'(z)>0$ and $R''(z)<0$ in the interior of $I_1$ by \cite{Br}. Then it is easy to see that $z_n=R^{-n}_1(z)$ converges to zero exponentially fast. Hence 
$$
cz_{n+1}=cR^{-1}_1(z_n)=z_n+O\big(z_n^2\big)_{n\to\infty}
=z_n\Big(1+O\big(z_n\big)_{n\to\infty}\Big)
$$
which proves (\ref{emathcalR}).

Let $W_n$ be the set of sequences (words) $w=k_1\dots k_n$ of length $n$, where $k_i\in\{1,\ldots,N\}$. 
For a word $w=k_1\ldots k_n\in W_n$ we write 
$R^{-n}_{w} = R^{-1}_{k_1}\dots R^{-1}_{k_n}$. 
We use the convention that 
$W_0=\{\varnothing\}$ and 
$R^{0}_{w}(z)=z$ for $w=\varnothing\in W_0$. 

Let 
$W_n^+$ be a subset of $W_n$ of words that do not start with $1$. 
Then we have 
\BEQ{expression}{\begin{aligned}
\zetaa{z_0}   {R}(s)\ps&
=
\lim\limits_{n\to\infty}
\sum_{z\in R^{-n}\{z_0\}}
(c^n z)^{- {s} /2}
\\
&
=\lim\limits_{n\to\infty} 
\sum_{m=0}^n c^{- {s} m/2} 
\sum_{w\in W_m^+} 
\left(c^{n-m} R^{-n+m}_1 R^{-m}_w(z_0)\right)^{- {s} /2}  
\\
&
=\sum_{m=0}^\infty c^{- {s} m/2} 
\sum_{w\in W_m^+} 
\left(\mathcal R R^{-m}_w(z_0)\right)^{- {s} /2} 
\end{aligned}}
for $\text{Re}(s)\ps>d_R$. 

Here and in what follows, in order to simplify notation, we use the convention that 
$z^{- {s} /2}=0$ for all real $s$ if $z=0$.

The last expression in (\ref{expression}) can be rewritten as 
\BEQ{expression2}{
\zetaa{z_0}   {R}(s)\ps=
\sum_{m=0}^\infty 
\sum_{w\in W_m^+} 
\left(c^{m} \mathcal R \left(R^{-m}_w(z_0)\right)\right)^{- {s} /2} 
=\sum_{j=1}^\infty
\lambda_j^{- {s} /2}
}
where the nondecreasing unbounded sequence of positive real numbers $\{\lambda_j\}_{j=1}^\infty$ 
is obtained by reordering the countable locally finite set of positive real numbers 
$
\bigcup_{m=0}^\infty 
\bigcup_{w\in W_m^+} 
\left\{c^{m} \mathcal R \left(R^{-m}_w(z_0)\right)\right\} \backslash\{0\}
$, taking the multiplicities into account.
\rP

\BTHM{thmPoly}{ 
The zeta function 
$\zetaa{z_0}   {R}(s)$ of the 
polynomial $R(z)$ has a meromorphic continuation of the 
form 
$$ 
\zetaa{z_0}   {R}(s)=\frac{\fone(s)}{1-Nc^{- {s} /2}}+\ftwo{z_0}(s)
$$
where $\fone(s)$ is analytic in $\mathbb C$ 
and $\ftwo{z_0}(s)$ is analytic for 
$ \text{Re}(s)\pms>0$ and is defined by \eqref{e-xiR}. 
The set of poles of $\zetaa{z_0}   {R}(s)$ is contained in 
$$ 
\left\{d_R+\tfrac{4in\pi}{\log c}\pms: n\in\mathbb Z\right\}
,$$ and there always is a pole at $d_R$. 

If $\mathcal J_R$ is totally disconnected, 
then there exists $\varepsilon\pms>0$ such that 
$\ftwo{z_0}(s)$ is analytic for 
$ \text{Re}(s)\pms> -\varepsilon$. 
If, in addition, $r=\max\limits_{z\in\mathcal J_R}| R'(z)|^{-1} <1$ then 
we can take 
$$ \varepsilon = -\frac{2\log r}{\log c}.$$}

Independently of our work this result was recently  improved in \cite{DGV}. In particular, it is shown 
in  \cite{DGV} that $\zetaa{z_0} {R}$  has a meromorphic continuation to the 
entire complex plain for any polynomial $R$ satisfying the same conditions as in our paper. 
However, an advantage of our method  is that it is 
 also applicable in the case when $R$ is a rational function satisfying 
appropriate assumptions. Therfore it is likely that our results can be extended to 
the class of fractals considered in \cite{MT}.

\Pr[Proof of Theorem~\ref{thmPoly}] 
To the last expression in (\ref{expression}) we add and subtract 
$$
\frac{\fone(s)}{1-Nc^{- {s} /2}}
$$
where
\BEQ{e-xiR}{
\fone(s)= 
\int_{\mathcal J_+} (\mathcal R (z))^{- {s} /2} d\kappa(z)
.} 
Here we use the fact (see, for instance \cite{BGH,BGM84,BGM88,BGM90,Br}) that on $\mathcal J_R$ there is a unique balanced invariant probability measure 
$\kappa$ which satisfies $$\kappa(R^{-1}_{j_1}\cdots R^{-1}_{j_m}(I))=\frac1{N^{m}}$$ 
for any sequence ${j_1},...,{j_m}$.

Note that $\fone(s)$ is an analytic function since $\mathcal J_+$ is 
separated from zero, and that 
$\frac{1}{1-Nc^{- {s} /2}}=\sum_{m=0}^\infty 
N^m c^{- {s} m/2} $. 
Then we only have to show that 
\begin{multline*}\label{eRem}
\sum_{m=0}^\infty 
c^{- {s} m/2} 
\left( 
\sum_{w\in W_m^+} 
\left(\mathcal R R^{-m}_w(z_0)\right)^{- {s} /2} 
- 
N^m \int_{\mathcal J_+} 
(\mathcal R (z))^{- {s} /2} d\kappa(z)
\right) 
=\\ 
\sum_{m=0}^\infty 
N^{m} c^{- {s} m/2} 
\sum_{w\in W_m^+} 
\int\limits_{\hskip5pt R^{-m}_w\{\mathcal J_+\}\hskip-5pt} 
\left[\left(\mathcal R R^{-m}_w(z_0)\right)^{- {s} /2} 
-
(\mathcal R (z))^{- {s} /2} \right] d\kappa(z)
\end{multline*}
is analytic when $ \text{Re}(s)\pms>0$. 

Thus we have to show that the series converges 
uniformly when $s$ is in a compact subset. 
To this end we can estimate 
the absolute value of the expression in square brackets by 
$$M 
\max_{{z\in R^{-m}_w\{\mathcal J_+\}}}
|z-R^{-m}_w(z_0)|\leqslant M \lambda(R^{-m}_w\{I_0\})$$ 
where $\lambda$ is the Lebesgue measure on $\mathbb R$, 
$$
M=
\max_{z\in I_+}
\left|\frac{d}{dz}\big(\mathcal R (z)\big)^{- {s} /2}\right|
,$$ 
and \,\,$I_+=\bigcup_{w\in W_{1}^+} 
{ R^{-1}_w\{I_0\}}$.

Then, since $\kappa(R^{-m}_w\{\mathcal J\})=N^{-m}$, 
\begin{multline*}
\sum_{w\in W_m^+} 
\int\limits_{\hskip5pt R^{-m}_w\{\mathcal J_+\}\hskip-5pt} 
\left|\left(\mathcal R R^{-m}_w(z_0)\right)^{- {s} /2} 
-
(\mathcal R (z))^{- {s} /2} \right| d\kappa(z) 
\leqslant \\ 
MN^{-m} \lambda(R^{-m+1}(I))
.\end{multline*} 
In this estimate we have to use the exponent $-m+1$ instead of $m$ 
because we assumed that $z_0\in I_0$, not 
necessarily $z_0\in I$. 

Thus we have the required convergence for 
$ \text{Re}(s)\pms>0$. 

If $\mathcal J_R$ is totally disconnected, 
then by pages 122-124 in \cite{Br}, there is $r<1$ such that 
$$
\lambda(R^{-n}(I))<const\cdot r^n.
$$ 
This estimate gives the required convergence for 
$ \text{Re}(s)\pms>-\varepsilon=\frac{2\log r}{\log c}$. 
Finally, note that $r\leqslant\max\limits_{z\in\mathcal J_R}| R'(z)|^{-1}$. 
\rP

In the applications of this theorem considered in this paper, 
the assumptions on the 
polynomial $R(z)$ given above are always satisfied because the Julia set of $R(z)$ is 
contained in the spectrum of a nonnegative self-adjoint operator (see \cite{MT,T01,T02}). 
Moreover, 
\begin{equation}\label{ePZF2}
\zetaa{z_0}   {R}(s)=
\sum\lambda_j^{-{ {s} }\sm{/2}}
\end{equation}
where $\{\lambda_j\}$ are eigenvalues of a Laplacian, and so 
the spectral zeta function has a representation as a \zf of a 
polynomial. In these cases the parameter $z_0$ usually depend on the boundary conditions. 
However, even if the polynomial $R(z)$ is not related to any Laplacian, 
the representation \eqref{ePZF2} still holds with 
 $\{\lambda_j\}$ defined in terms of 
the limit points of $\bigcup_{n>0}c^nR^{-n}(z_0)$ as in \Thm{propPZF}.

\BTHM{thmRiemann}{The Riemann zeta function $\zeta(s)$ has a representation 
$$
\zeta(s)= 
\frac12C^s
\zetaa{0}   {R}(s)
$$ 
where $C=\sqrt2\pi$ and $\zetaa{0}   {R}(s)$ is 
the \zfotp 
$R(z)=2z(2-z)$.}

\Pr For the Neumann or Dirichlet Laplacian $\Delta$ on $I=[0,1]$ we have 
$$
\zeta(s)=
\sum_{j=1}^\infty{j}^{-s}=
\pi^{ s} \sum_{j=1}^\infty\left({\pi^2}{j^2}\right)^{-s/2}=
\pi^{ s} \zeta_\Delta(s).
$$
Let 
$$
\Delta_n f\big(\tfrac{k}{2^n}\big)=
f\big(\tfrac{k}{2^n}\big)-\tfrac12 f\big(\tfrac{k-1}{2^n}\big)
-\tfrac12 f\big(\tfrac{k+1}{2^n}\big)
$$
be the probabilistic Laplacian on $\{0,\ldots,\tfrac{k}{2^n},\tfrac{k+1}{2^n},\ldots,1\}$, which is  
 the generator of the simple nearest neighbor random walk. 
Then 
$$\Delta f(x)=-f''(x)=
2
\lim\limits_{n\to\infty}4^n\Delta_nf(x)
$$ 
for 
any twice continuously differentiable function $f$. 
The theorem then follows from the definition of 
$\zetaa{0}   {R}(s)$ and 
Lemma~\ref{lem06} since 
$$\zeta_\Delta(s) =\frac122^{-s/2}4^{s/2}\zetaa{0}   {R}(s)=\frac122^{s/2}\zetaa{0}   {R}(s).$$ 
 Note that the preimages of zero are $\{0,2\}=\sigma(\Delta_\ii {0})$, and so the spectrum of~$\Delta_\ii {n}$ 
as a set coincides with $R^{-n-1}(0)$, which explains the extra factor $4^{s/2}$ above.   Also note 
that $1$ is the preimage of $2$ of multiplicity two, but $1$ as a simple  eigenvalue of~$\Delta_\ii {n}$, which explains the extra factor $\frac12$ above.  
\rP

\BLEM{lem06}{
If $z\neq1$ then 
$ R(z)=2z(2\ms-z)\in\sigma(\Delta_\ii {n})$ if and only if $
z\in\sigma(\Delta_\ii {n\ps+1})$, with the same multiplicities. Moreover,  the Neumann discrete 
\Lp s have  simple spectrum with 
$\sigma(\Delta_\ii {0})=\{0,2\}$ and $1\in\sigma(\Delta_\ii {n})$ for all $n>0$.}
\Pr This lemma can be proved by an elementary computation, see   \cite[Example 2.4 and Theorem 5.8]{MT}. 
\rP

\section{Mellin transform}\LL{sMellin}

In this expository section we contrast the \zfoap\ with the 
Mellin transform of an invariant measure related to the same polynomial. 
We do not use the correspondence between 
these object in our 
work, but the reader can find some related formulas  in \cite{DGV,Grabner}. 

\BDF{dfMellin}{
If $\kappa$ is the unique balanced invariant measure on 
the Julia set $\mathcal J_R$ of $R(z)$ then its 
\bydef{Mellin transform} is
$$
\M\ind{R}(s)\ps=
\lim\limits_{n\to\infty}N^{-n}
\lala\lala\lala\lala\lala\lala\lala
\sum_{z\in R^{-n}\{z_0\}}
\lala\lala\lala\lala
 z^{- {s} /2}
\ps=\int_{\mathcal J_R} 
\lala\lala\lala
z^{- {s} /2}d\kappa(z).
$$
}

\BBF{figJR}
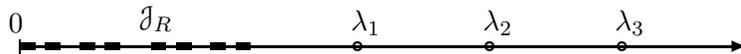
{A pictorial illustration of the Julia set ${\mathcal J_R}$ of $R(z)$ and the sequence $\{\lambda_j\}$ in \Thm{propPZF}.}
{\begin{picture}(275,25)(0,0) \setlength{\unitlength}{0.5pt}\small
\thicklines
\put(0,0){\vector(1,0){550}}
\put(-9,11){${0}$}
\put(90,15){${\mathcal J_R}$}
\linethickness{5.5pt}%
\put(-1,0){\line(1,0){1}}
\put(0,0){%
\linethickness{2.5pt}%
\put(0,0){\line(1,0){11}}
\put(19,0){\line(1,0){11}}
}%
\put(45,0){%
\linethickness{2.5pt}%
\put(0,0){\line(1,0){11}}
\put(19,0){\line(1,0){11}}
}%
\put(99,0){
\put(0,0){%
\linethickness{2.5pt}%
\put(0,0){\line(1,0){11}}
\put(19,0){\line(1,0){11}}
}%
\put(45,0){%
\linethickness{2.5pt}%
\put(0,0){\line(1,0){11}}
\put(19,0){\line(1,0){11}}
}%
}%
\put(250,12){${\lambda_1}$}
\linethickness{5.5pt}%
\put(255,0){\circle{5}}
\put(350,12){${\lambda_2}$}
\linethickness{5.5pt}%
\put(355,0){\circle{5}}
\put(450,12){${\lambda_3}$}
\linethickness{5.5pt}%
\put(455,0){\circle{5}}
\end{picture}}


\BBF{figZFMT}
{A pictorial illustration of the poles of the \zfoap\ in \Thm{thmPoly} (left), and of the poles of its Mellin transform (right).}
{\begin{picture}(120,170)(0,-90) \setlength{\unitlength}{0.5pt}\small\scriptsize
\thicklines\linethickness{.99pt}
\setlength{\unitlength}{0.4pt}
\put(205,7){$\sm{0}$}
\put(205,5){\line(0,-1){10}}
\put(-0,0){\vector(1,0){330}}
\put(75,7){$-\lg\varepsilon$}
\put(275,7){$\sm{d_R}$}
\multiput(270,-180)(0,36){11}{\circle{7}}
\multiput(70,-200)(0,20){20}{\line(0,1){10}}
\thinlines
\multiput(70,-150)(0,4){87}{\line(-5,-3){90}}
\end{picture} 
\qquad 
\qquad 
\begin{picture}(140,170)(0,-90) \setlength{\unitlength}{0.5pt}\small\scriptsize
\thicklines\linethickness{.99pt}
\setlength{\unitlength}{0.4pt}
\put(20,7){$\sm{0}$}
\put(20,5){\line(0,-1){10}}
\put(0,0){\vector(1,0){370}}
\put(95,7){$\sm{d_R}$}
\multiput(90,-180)(0,36){11}{\circle{7}}
\put(195,7){$\sm{d_R\ps+1}$}
\multiput(190,-180)(0,36){11}{\circle{7}}
\put(295,7){$\sm{d_R\ps+2}$}
\multiput(290,-180)(0,36){11}{\circle{7}}
\end{picture}}

 \BTHMm{\cite{BGM84}}{
The Mellin transform $\M\ind{R}$ converges for $\text{Re}(s)\ps<d_R$ and 
has a meromorphic continuation to $\mathbb C$ 
with the set of poles contained in 
$$ 
\left\{\frac{2\log N+4in\pi}{\log c}
\ps+m\ps: n\in\mathbb Z,\ m\in\mathbb N\right\}
.$$
Every real pole has a positive residue and thus is not 
degenerate. 
}

Note that by the Theorem of Bessis, Geronimo and Moussa 
the Mellin transform has infinitely many nondegenerate 
positive real poles, 
while by Theorem~\ref{thmPoly} 
the zeta function of a polynomial has only one 
positive real pole (see Figure~\ref{figZFMT}).

\NEWSECTION3{Fractal Laplacian on an interval and its zeta function}{s4}
\def\x{\sm{(x)}}
\def\ih#1{{{\huge\mbox{$#1$}}}}
\def\ih#1{$#1$}

We begin with defining the usual energy form of the unit interval 
by introducing a self-similar structure. In this way the interval can be 
considered as a p.c.f. fractal (see \cite{Ba,BNT,Ki1,Ki}). 
In particular, in these papers the reader can find proofs of the propositions given below, 
as well as the general theory of Dirichlet forms on fractals and further references. 

We use three contractions because it will allow us to easily modify the construction in order to obtain a fractal \Lp\ on the unit interval. One can make the same construction using  any number of contractions. Note that in Theorem 
\ref{thmRiemann} we essentially used a self-similar structure on the unit interval with two contractions. 

 \begin{remark}
We emphasize that in this and the next section we consider two different 
iteration function systems which are not to be confused. We start with  an iteration 
function system of linear contractions to define a self-similar set, which is the set where the \Lp\ is defined (and the state space of the corresponding diffusion process). 
This self-similar set is often a fractal, but in this section it is an interval with a self-similar structure. 
Then, based on the ``spectral decimation'' or, more appropriately called, ``spectral self-similarity'' (see \cite{MT}), we describe the spectrum of the \Lp\ using inverse iterations of a polynomial (a rational function in more general situation of \cite{MT}). To this latter iteration function system we apply the results of Section~\ref{s3}. 
\end{remark}

We define three contractions 
$\,F_1,F_2,F_3:\mathbb R^1\to\mathbb R^1\,$ 
with fixed points $x_j=0,\frac12,1$ by 
 $$F_j(x)=\tfrac13x+\tfrac23x_j.$$ 
Then the interval $I{=}[0,1]$ 
is a unique compact 
 set such that 
$$I=
\bigcup\limits_{\sm{\sm{{j{=}1,2,3}}}}\hskip-.975ex
F_j(I).$$
The \bydef{discrete approximations} to $I$ are 
defined inductively by 
$$V_n = 
\hskip-.975ex\bigcup\limits_{\sm{\sm{{j{=}1,2,3}}}}\hskip-.975ex
F_j(V_{n\sm{\sm{-}}1})
=\big\{\tfrac k{3^n}\big\}_{k{=}0}^{3^n}$$ 
where $V_0=\partial I=\{0,1\}$ is the
\bydef{boundary} of $I$. For $x,y\in V_\ii {n}$ we write $y \sim x$ if $|x-y|=3^{-n}$. 
 

\BDF{df07}{The \bydef{discrete Dirichlet (energy) form} on $V_n$ is 
 $$
\mathcal E_n (\ddD{f}) = 3^ n
{\sum 
}
\big(f\sm{(y)} {-} f\x \big)^\ii 2
$$ 
where the sum is taken over all unordered pairs $x,y\in V_\ii {n}$ 
such that $y \sim x$. 
The \bydef{Dirichlet (energy) form} on $I$ is 
$$
\mathcal E (\ddD{f}) = \lim\limits_\ii {n{\to}\infty} \mathcal E_n (\ddD{f}) 
$$
if this limit exist and is finite. 
By the next proposition this limit is increasing and so it always exists as a nonnegative real number or $+\infty$. 

A function $h$ is \bydef{harmonic} if it minimizes the 
energy given the boundary values.}

\BPROP{prop09}{We have that 
$
\mathcal E_{n{+}1} (\ddD{f})\geqslant\mathcal E_n (\ddD{f})
$ 
 for any function $f$, and 
$$
\mathcal E_{n{+}1}(\ddD{h})= 
\mathcal E_n (\ddD{h})=\mathcal E(\ddD{h})
$$ 
for a harmonic $h$. A function $h$ is harmonic if and only if it is linear. If $f$ is continuously differentiable 
then $$
\mathcal E (\ddD{f}) =\int_0^1|f'(x)|^2 dx.
$$
The Dirichlet (energy) form $\mathcal E$ on $I$ is \bydef{self-similar} 
in the sense that 
$$ 
\mathcal E(\ddD{f})=3
\sum_{\sm{\sm{j=1,2,3}}}
\mathcal E(\ddD{f\sm{\circ}F_j}).
$$} 

\BDF{df11}{The \bydef{discrete \Lp s} on $V_n$ are 
$$ 
\Delta_n f\x = \tfrac12 
{\sum\limits_ {\genfrac{}{}{0pt}{2}{y\in V_\ii {n}}{y \sim x}}}
 \big(f\sm{(y)} {-} f\x \big), \quad x{\in}V_n\backslash V_0
$$ 
and the (renormalized) \Lp\ on $I$ is 
$$
\Delta f\x = \lim\limits_\ii {n{\to}\infty} 9^ n 
\Delta_n f\x =\frac12f''\x 
$$ 
for any twice differentiable function.}
Obviously,  for any twice differentiable function  this \Lp\ satisfies the \bydef{Gauss--Green (integration by parts) formula}
$$
 \mathcal E (\ddD{f})=-2\int_0^1f\Delta f dx + ff'\Big|^1_0 
$$
and has the \bydef{spectral asymptotics}
$$
\lim_{\lambda{\to}\infty}\frac{\rho(\lambda)}{\lambda^{d_s/2}}=\frac1\pi
$$
where $\rho(\lambda)$ is the \bydef{eigenvalue counting function} 
$
\rho(\lambda)=\#\{j:\lambda_j<\lambda\}
$
of the Dirichlet or Neumann \Lp\ $\Delta$ 
and $d_s=1$ is its spectral dimension.

In what follows we modify this construction to obtain a one parameter family 
of fractal \Lp s on the unit interval. We introduce a parameter $0<p<1$ and define $q=1-p$, which later will be shown to have a meaning of transitional probabilities of a random walk. 
Then we define contraction factors (or ``resistance weights'') 
$$ 
r_1=r_3=\frac{p+pq}{2+pq}\text{ \ \ and \ \ } 
r_2=\frac{2q-pq}{2+pq}, 
$$
and 
``measure weights''
$$ 
m_1=m_3=\frac{q}{1+q}\text{ \ \ and \ \ } 
m_2=\frac{p}{1+q}. 
$$
Note that, up to a constant factor, the resistance weights are reciprocals 
of the measure weights, and 
$$m_1+m_2+m_3=r_1+r_2+r_3=1.$$ By their definition, 
constants resistance weights $r_i$ and  measure weights $m_i$ 
depend on the choice of $p$.

We define three contractions 
$\,F_{p,1},F_{p,2},F_{p,3}:\mathbb R^1\to\mathbb R^1\,$ 
with fixed points $x_j=0,\frac12,1$ by 
 $$F_{p,j}(x)=r_jx+(1-r_j)x_j.$$ 
Then the interval $I{=}[0,1]$ 
is a unique compact 
 set such that 
$$I=
\bigcup\limits_{\sm{\sm{{j{=}1,2,3}}}}\hskip-.975ex
F_{p,j}(I).$$
The \bydef{discrete approximations} to $I$ are 
defined inductively by 
$$V_{p,n} = 
\hskip-.975ex\bigcup\limits_{\sm{\sm{{j{=}1,2,3}}}}\hskip-.975ex
F_{p,j}(V_{p,n\sm{\sm{-}}1})
$$ 
where, as before, $V_0=\partial I=\{0,1\}$ is the
\bydef{boundary} of $I$. 

\BDF{df07-}{The \bydef{discrete Dirichlet (energy) form} on $V_n$ is 
defined inductively 
$$ 
\mathcal E_{p,n} (\ddD{f})=
\sum_{\sm{\sm{j=1,2,3}}}
\tfrac1{r_j} \mathcal E_{p,n-1}(f\sm{\circ}F_j).
$$
with $
\mathcal E_0 (\ddD{f}) = 
\big(f(1)-f(0)\big)^2 $ as before.  
The \bydef{Dirichlet (energy) form} on $I$ is 
$$
\mathcal E_p (\ddD{f}) = \lim\limits_\ii {n{\to}\infty} 
\mathcal E_{p,n} (\ddD{f}) 
$$
if this limit exist and is finite. Note that this is the same definition as the Dirichlet form defined before 
if and only if $p=\frac12$.
By the next proposition this limit is increasing and so it always exists as a nonnegative real number or $+\infty$.

The next two propositions follow, for example, from the general results in \cite{Ba,BNT,Ki1,Ki}. Again, a function $h$ is \bydef{harmonic} if it minimizes the 
energy given the three boundary values.}

\BPROP{prop099}{We have that 
$
\mathcal E_{p,n{+}1} (\ddD{f})\geqslant\mathcal E_{p,n} (\ddD{f})
$ 
 for any function $f$, and 
$$
\mathcal E_{p,n{+}1}(\ddD{h})= 
\mathcal E_{p,n} (\ddD{h})=\mathcal E_p(\ddD{h})
$$ 
for a harmonic $h$.}

\BPROP{prop16}{
The Dirichlet (energy) form $\mathcal E_p$ on $I$ is local and regular, and is \bydef{self-similar} 
in the sense that 
$$ 
\mathcal E_p(\ddD{f})=
\sum_{\sm{\sm{j=1,2,3}}}
\tfrac1{r_j} \mathcal E_p(\ddD{f\sm{\circ}F_{p,j}}).
$$
The domains of $\mathcal E_p$ and of the corresponding \Lp\ $\Delta_\mu $ 
are contained and dense in the space of continuous functions. 

There exists a \bydef{$\mu$--\Lp} $\Delta_\mu$, that is a densely defined 
self-adjoint operator that satisfies  
the following 
Gauss--Green (integration by parts) formula 
$$
 \mathcal E_p (\ddD{f})=
{-}C\int_0^1f\Delta_\mu f d\mu + ff'\big|^1_0. 
$$
where $\mu$ is a unique probability \bydef{self-similar} 
measure 
with weights $m_1,$ $m_2,$ $m_3$, that is 
$$ 
\mu =\sum_{\sm{\sm{j=1,2,3}}}m_j\, \mu \sm{\circ} F_{p,j}.
$$
Also
$$
\Delta_\mu f\x 
= 
\hskip-.21em
\lim\limits_\ii {n{\to}\infty} 
\hskip-.21em
\big(1\sm{+}\tfrac2{pq}\big)^ n 
\Delta_{p,n} f\x 
$$
where the discrete \Lp s 
$$
\Delta_{p,n} f(y_k)=\left\{
\begin{aligned}
\mbox{$pf(y_{k-1})
+qf(y_{k+1}) -f(y_k)$}&\\ 
& \ \ \text{or} \\
\mbox{$qf(y_{k-1})
+pf(y_{k+1}) -f(y_k)$}& \end{aligned}
\right. \ 
$$
 are defined as the 
generators of the nearest neighbor random walks on $V_{p,n}$ 
with transitional 
probabilities $p$ and $q$ assigned according to the 
weights of the corresponding intervals.}


Note that by definition 
$p=\frac{m_2}{m_1{+}m_2}$, \,$q=\frac{m_1}{m_1{+}m_2}$. 
The transitional 
probabilities $p$ and $q$ can be assigned inductively 
as shown in \Fig{figRW}. 

\BBF{figRW}
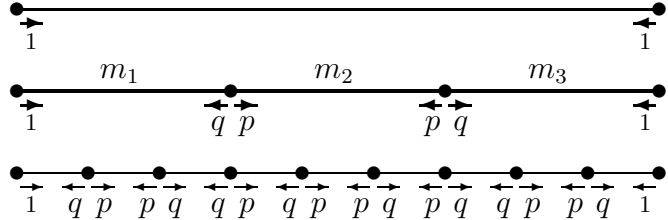
{Random walks corresponding to the discrete \Lp s $\Delta_{p,n}$.}
{\begin{center}
\begin{picture}(246,30)(0,-20) \setlength{\unitlength}{.45pt}
\thicklines
\put(0,0){\circle*{12}}
\put(540,0){\circle*{12}}
\put(0,0){\line(1,0){540}}
\put(3,-12){\vector(1,0){19}}
\put(537,-12){\vector(-1,0){19}}
\put(7,-33){{\scriptsize\ih1}}
\put(523,-33){{\scriptsize\ih1}}
\end{picture}

\begin{picture}(246,30)(0,-20) \setlength{\unitlength}{.45pt}
\thicklines
\put(0,0){\circle*{12}}
\multiput(0,0)(60,0){9}{\line(1,0){60}}
\multiput(180,0)(180,0){3}{\circle*{12}}
\multiput(3,-12)(180,0){3}{\vector(1,0){19}}
\multiput(177,-12)(180,0){3}{\vector(-1,0){19}}
\put(7,-33){{\scriptsize\ih1}}
\put(70,11){\ih{m_1}}
\put(250,11){\ih{m_2}}
\put(430,11){\ih{m_3}}
\put(163,-33){\ih q}
\put(187,-33){\ih p}
\put(343,-33){\ih p}
\put(367,-33){\ih q}
\put(523,-33){{\scriptsize\ih1}}
\end{picture}

\begin{picture}(246,30)(0,-20) \setlength{\unitlength}{.45pt}%
\put(0,0){\circle*{12}}
\multiput(0,0)(60,0){9}{\line(1,0){60}}
\multiput(60,0)(60,0){9}{\circle*{12}}
\multiput(3,-12)(60,0){9}{\vector(1,0){19}}
\multiput(57,-12)(60,0){9}{\vector(-1,0){19}}
\put(7,-33){{\scriptsize\ih1}}
\put(43,-33){\ih q}
\put(67,-33){\ih p}
\put(103,-33){\ih p}
\put(127,-33){\ih q}
\put(163,-33){\ih q}
\put(187,-33){\ih p}
\put(223,-33){\ih q}
\put(247,-33){\ih p}
\put(283,-33){\ih p}
\put(307,-33){\ih q}
\put(343,-33){\ih p}
\put(367,-33){\ih q}
\put(403,-33){\ih q}
\put(427,-33){\ih p}
\put(463,-33){\ih p}
\put(487,-33){\ih q}
\put(523,-33){{\scriptsize\ih1}}
\end{picture}
\end{center}
}

Note that the previous construction of the standard \Lp\ and the 
Dirichlet form on $I=[0,1]$ 
corresponds to 
$p=\frac12$. If $p\not=\frac12$ then 
we can make a change of variable on the unit interval $I$ so that 
either the Dirichlet form becomes the standard one, or the measure 
becomes the Lebesgue measure, but not both. Thus for different values of 
$p$ the $\mu$--\Lp s are different even up to a change of variable.

By the result of Kigami and Lapidus \cite{KL}, 
the Dirichlet or Neumann \Lp\ $\Delta_\mu$ 
has the \bydef{spectral asymptotics}
$$
0<
\liminf_{\lambda{\to}\infty}
\frac{\rho_p(\lambda)}{\lambda^{d_s/2}}
\leqslant
\limsup_{\lambda{\to}\infty}
\frac{\rho_p(\lambda)}{\lambda^{d_s/2}}
<\infty
$$
where $\rho_p(\lambda)$ is the eigenvalue counting function of $\Delta_\mu$, and its spectral dimension is 
$$
d_s=\dfrac{\log9}{\log\big(1\sm{+}\tfrac2{pq}\big)}\leqslant1,
$$ 
where the inequality is strict if and only if $p\neq q$.


In what follows we denote  
$$R_p(z)\ps=z( z^2 \pms+3z\pms+2\pms+pq)/pq.$$ 
We use notation $\sigma(A)$ for the spectrum of a linear operator $A$.

\BLEM{prop17}{If $z\neq-1\pm p$ 
then $ R(z)\in\sigma(\Delta_\ii {p,n})$ 
if and only if 
$z\in\sigma(\Delta_\ii {p,n\ps+1})$, with the same multiplicities. 
Moreover,  the Neumann discrete 
\Lp s have  simple spectrum with 
$\sigma(\Delta_\ii {0})=\{0,-2\}$ 
and
$$
\sigma(\Delta_\ii {p,n})=\{0,-2\}\bigcup_{m=0}^{n-1}R_p^{-m}\{-1\pm q\}
$$
for 
all $n>0$. 
In particular, for all $n>0$ we have  $0,-1\pm q,-2\in\sigma(\Delta_\ii {p,n})$. 
Also, for 
all $n>0$ we have  $-1\pm p\in\sigma(\Delta_\ii {p,n})$ 
if and only if $p=q$.} 

\Pr 
In this case, according to \cite[Lemma 3.4]{T}, \cite[(3.2)]{MT} or  \cite{T01,T02}, we have that $R_p(z)=\dfrac{\varphi_1(z)}{\varphi_0(z)}$ where 
$\varphi_0$ and $\varphi_1$ solve the matrix equation 
$$
S - z I_0 -  \bar   X    (Q-z I_1)^{-1}   X    =\varphi_0 (z) H_0  
-\varphi_1 (z) I_0.
$$
with $S=\bar X=X=I_0=I_1$, the $2\times2$ identity matrices, 
$$
Q=\left(\begin{array}{rr}
 -1 &  p\cr
 p & -1\cr
\end{array} \right)
$$
and
$$
H_0=\left(\begin{array}{rr}
 -1 &  1\cr
 1 & -1\cr
\end{array} \right).
$$
Solving this we obtain 
$$
\varphi_0(z)=\frac{pq}{z^2+2z+1-p^2}
$$
and 
$$
\varphi_1(z)=\frac{z( z^2 \pms+3z\pms+2\pms+pq)}{z^2+2z+1-p^2}
.$$
 Then we use the abstract spectral self-similarity result (see \cite{T,MT,T01,T02}) that 
the spectrum of $\Delta_\ii {p,n+1}$ can be obtain from 
the spectrum of $\Delta_\ii {p,n}$ by the inverse  of  $R(z)$. 
Note that $0$ and $2$ are fixed points of $R(z)$. 
The preimages of $0$ are $0$, $-1-p$ and $-1-q$. 
The preimages of  $2$ are $2$, $-1+p$ and $-1+q$. If $p\neq q$ then $-1\pm p$ are not 
eigenvalues because they are poles of $\varphi_0(z)$ (see \cite{MT,T,T01,T02}). 
\rP

\BTHM{thm18}{
If $p\neq \frac12$ then $d_s=d_R<1$ and the \szf\ of $\Delta_\mu$ is 
$$\zeta_{\Delta_\mu}(s)=
\frac{c^ {s} }{ 1-c^{-s/2} }\left(\zetaa{z_1}   {R}(s)+\zetaa{z_2}   {R}(s)\right)$$ 
where $c=1\sm{+}\tfrac2{pq}$ and $z_1,z_2=1\pm q$. 
Moreover, $\zeta_{\Delta_\mu}(s)$ has a 
meromorphic continuation to 
$ \text{Re}(s)>-\varepsilon$ for some  $\varepsilon_p>0$, and  
 $\liminf\limits_{p\to0}\varepsilon_p\geqslant2$, 
$\liminf\limits_{p\to1}\varepsilon_p\geqslant2$.} 

\Pr 
This follows from  \Lem{prop17}, \Thm{thmPoly} and the definition of the \zfoap. 
The last statement follows from the 
fact that as $p\to0$ or $p\to1$ we have $\lim\min\limits_{z\in \mathcal J_R}|R_p'(z)|=\infty$. 
 \rP

\BTHM{thmRiemann3}{The Riemann zeta function $\zeta(s)$ 
has a representation 
$$
\zeta(s)=
\frac12C^s \zetaa{0}   {R}(s)
$$
where $C=\frac{3\pi}{\sqrt2}$ and $\zetaa{0}   {R}(s)$ is 
the \zfotp 
 $$R(z)= z( 4z^2 -12z+9). $$} 

\Pr 
The proof follows from  \Lem{prop17}, \Thm{thmPoly} and the definition of the \zfoap,  
 analogously to the proof of Theorem~\ref{thmRiemann}, and of Theorem~\ref{thm18} with $p=q$. 
\rP

\BBF{figCD}
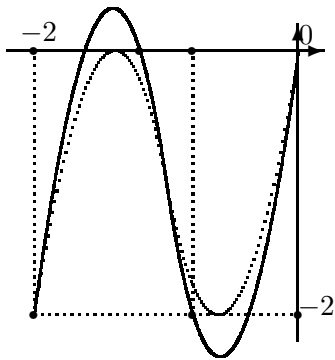
{Sketch of the cubic polynomial $R_p(z)$ 
associated with the fractal \Lp s on the interval.}
{\begin{picture}(120,133)(-50,-66) \setlength{\unitlength}{0.5pt}\small
\thicklines
\put(-120,100){\vector(1,0){240}}
\put(100,-120){\vector(0,1){240}}
\put(101,105){$\ii{0}$}
\put(-110,107){$\ii{-2}$}
\put(100,-100){$\ii{-2}$}
\put(-100,100){\circle*{5}}
\put(-20,100){\circle*{5}}
\put(20,100){\circle*{5}}
\put(20,-100){\circle*{5}}
\put(100,-100){\circle*{5}}
\put(-100,-100){\circle*{5}}
\qbezier[35](-100,-100)(-100,0)(-100,100)
\qbezier[35](20,-100)(20,0)(20,100)
\qbezier[35](-100,-100)(0,-100)(100,-100)
\qbezier(-100,-100)(-45,307)(0,0)
\qbezier[63](-100,-100)(-45,241.5)(0,0)
\qbezier[35](100,100)(0,100)(-100,100)
\qbezier(100,100)(45,-307)(0,0)
\qbezier[63](100,100)(45,-241.5)(0,0)
\end{picture}}

\BREM{remFI}
{In \Fig{figCD} we give a sketch that describes 
the complex dynamics of the family of cubic polynomials 
associated with the fractal \Lp s on the interval. 
Note that the Julia set is real, and so what we sketch is actually the real dynamics 
of these polynomials. 
The curved dotted line corresponds to the case when $p=\frac12$ and 
the Julia set of $R(z)=R_{\frac12}(z)$ is the interval $[-2,0]$. For any other value of $p$ ($0<p<1$), the graph of the polynomial $R_p(z)$ behaves like the shown solid curved line. It is easy to see that then the 
Julia set of $R_p(z)$ is a Cantor set of Lebesgue measure zero. 
Note that the change $p\mapsto1-p$ does not change the polynomial 
$R_p(z)$, although the \Lp s ${\Delta_\mu}$ are different.} 


\NEWSECTION4{Spectral zeta function of the Sierpi\'nski gasket}{s5} 

We begin this section with describing the analysis on the \Sig\ 
in the way 
first developed by J.~Kigami \cite{Ki0}. We recall only basic facts here, and 
in the papers \cite{BP,BGH,BST,DSV,DABK,Fu,FS,Ki1,Ki,Ko,MS,R,RT,T} 
the reader can find proofs of the propositions given below, as well as 
 the general theory of Dirichlet forms and diffusions on fractals, many physical applications 
and further references.

If we fix three contractions 
$F_\ii j(x)=\frac{1}{2}(x\pms+v_\ii j)$, 
where $v_j$ are vertices of a triangle in $\mathbb R^2$, 
then the \bydef{\Sig}
is a unique compact set $S$ 
such that 
$$S = F_\ii 1(S)\cup F_\ii 2(S)\cup F_\ii 3(S).$$

\def\LiNe#1#2#3#4#5{{
\count41=#3 \advance\count41 by -#1 
\count42=#4 \advance\count42 by -#2 
\count51=\count41 \ifnum\count41<0 \multiply\count41 by -1 \fi
\count52=\count42 \ifnum\count42<0 \multiply\count42 by -1 \fi 
\advance\count41 by \count42 \divide\count41 by \dvv 
\advance\count41 by 1 
\count42=\count41 
\count43=0 
\loop 
\count101=#1 
\count102=#2 
\count103=#3 
\count104=#4 
\multiply\count101 by \count42 \multiply\count103 by \count43 
\advance\count101 by \count103
\divide\count101 by \count41 
\multiply\count102 by \count42 \multiply\count104 by \count43 
\advance\count102 by \count104
\divide\count102 by \count41 
\put(\count101,\count102){#5}
\advance\count42 by -1 
\advance\count43 by 1 
\ifnum \count42>-1 \repeat
}}

\def\trigg#1#2#3#4#5#6#7{{
\count107=#1 
\count108=#2 
\advance\count107 by #3 \advance\count107 by #5 \divide\count107 by 3 
\advance\count108 by #4 \advance\count108 by #6 \divide\count108 by 3
\LiNe{#1}{#2}{\count107}{\count108}{#7} 
\LiNe{#3}{#4}{\count107}{\count108}{#7} 
\LiNe{#5}{#6}{\count107}{\count108}{#7} 
}}

\def\power#1#2{\count91=1 \count92=#2 \loop 
\advance\count92 by -1 \multiply\count91 by #1 
\ifnum \count92>0 \repeat}

\def\troo#1{\divide\count101 by 2 \divide\count102 by 2 
\ifcase#1 
 \or 
\advance\count101 by 8192 \or 
\advance\count101 by 4096 \advance\count102 by 7168 \fi} 

\def\SIGA#1#2{\power{3}{#1}
\loop 
\advance\count91 by -1 
\count92=#1 \count95=\count91 
\count101=0 \count102=0 
{\loop 
\advance\count92 by -1 
\count93=\count95 
\divide\count95 by 3 
\count94=\count95 
\multiply\count94 by 3 
\advance\count93 by -\count94 
\troo{\count93} 
\ifnum \count92>0 \repeat \put(\count101,\count102){#2}}
\ifnum \count91>0 \repeat} 

\def\trooh#1{
\ifcase#1 
%
%
%
\count141=0\count142=0\count143=0\count151=0\count152=0\count153=0
\count231=5 
\count232=0 
\count233=0 
\multiply\count231 by \count211 
\multiply\count232 by \count212 
\multiply\count233 by \count213 
\advance\count141 by \count231
\advance\count141 by \count232
\advance\count141 by \count233
\divide\count141 by 5 
\count231=5 
\count232=0 
\count233=0 
\multiply\count231 by \count221 
\multiply\count232 by \count222 
\multiply\count233 by \count223 
\advance\count151 by \count231
\advance\count151 by \count232
\advance\count151 by \count233
\divide\count151 by 5
\count231=2 
\count232=2 
\count233=1 
\multiply\count231 by \count211 
\multiply\count232 by \count212 
\multiply\count233 by \count213 
\advance\count142 by \count231
\advance\count142 by \count232
\advance\count142 by \count233
\divide\count142 by 5 
\count231=2 
\count232=2 
\count233=1 
\multiply\count231 by \count221 
\multiply\count232 by \count222 
\multiply\count233 by \count223 
\advance\count152 by \count231
\advance\count152 by \count232
\advance\count152 by \count233
\divide\count152 by 5
\count231=2 
\count232=1 
\count233=2 
\multiply\count231 by \count211 
\multiply\count232 by \count212 
\multiply\count233 by \count213 
\advance\count143 by \count231
\advance\count143 by \count232
\advance\count143 by \count233
\divide\count143 by 5 
\count231=2 
\count232=1 
\count233=2 
\multiply\count231 by \count221 
\multiply\count232 by \count222 
\multiply\count233 by \count223 
\advance\count153 by \count231
\advance\count153 by \count232
\advance\count153 by \count233
\divide\count153 by 5
\or
\count141=0\count142=0\count143=0\count151=0\count152=0\count153=0
\count231=2 
\count232=2 
\count233=1 
\multiply\count231 by \count211 
\multiply\count232 by \count212 
\multiply\count233 by \count213 
\advance\count141 by \count231
\advance\count141 by \count232
\advance\count141 by \count233
\divide\count141 by 5 
\count231=2 
\count232=2 
\count233=1 
\multiply\count231 by \count221 
\multiply\count232 by \count222 
\multiply\count233 by \count223 
\advance\count151 by \count231
\advance\count151 by \count232
\advance\count151 by \count233
\divide\count151 by 5
\count231=0 
\count232=5 
\count233=0 
\multiply\count231 by \count211 
\multiply\count232 by \count212 
\multiply\count233 by \count213 
\advance\count142 by \count231
\advance\count142 by \count232
\advance\count142 by \count233
\divide\count142 by 5 
\count231=0 
\count232=5 
\count233=0 
\multiply\count231 by \count221 
\multiply\count232 by \count222 
\multiply\count233 by \count223 
\advance\count152 by \count231
\advance\count152 by \count232
\advance\count152 by \count233
\divide\count152 by 5
\count231=1 
\count232=2 
\count233=2 
\multiply\count231 by \count211 
\multiply\count232 by \count212 
\multiply\count233 by \count213 
\advance\count143 by \count231
\advance\count143 by \count232
\advance\count143 by \count233
\divide\count143 by 5 
\count231=1 
\count232=2 
\count233=2 
\multiply\count231 by \count221 
\multiply\count232 by \count222 
\multiply\count233 by \count223 
\advance\count153 by \count231
\advance\count153 by \count232
\advance\count153 by \count233
\divide\count153 by 5
\or
\count141=0\count142=0\count143=0\count151=0\count152=0\count153=0
\count231=2 
\count232=1 
\count233=2 
\multiply\count231 by \count211 
\multiply\count232 by \count212 
\multiply\count233 by \count213 
\advance\count141 by \count231
\advance\count141 by \count232
\advance\count141 by \count233
\divide\count141 by 5 
\count231=2 
\count232=1 
\count233=2 
\multiply\count231 by \count221 
\multiply\count232 by \count222 
\multiply\count233 by \count223 
\advance\count151 by \count231
\advance\count151 by \count232
\advance\count151 by \count233
\divide\count151 by 5
\count231=1 
\count232=2 
\count233=2 
\multiply\count231 by \count211 
\multiply\count232 by \count212 
\multiply\count233 by \count213 
\advance\count142 by \count231
\advance\count142 by \count232
\advance\count142 by \count233
\divide\count142 by 5 
\count231=1 
\count232=2 
\count233=2 
\multiply\count231 by \count221 
\multiply\count232 by \count222 
\multiply\count233 by \count223 
\advance\count152 by \count231
\advance\count152 by \count232
\advance\count152 by \count233
\divide\count152 by 5
\count231=0 
\count232=0 
\count233=5 
\multiply\count231 by \count211 
\multiply\count232 by \count212 
\multiply\count233 by \count213 
\advance\count143 by \count231
\advance\count143 by \count232
\advance\count143 by \count233
\divide\count143 by 5 
\count231=0 
\count232=0 
\count233=5 
\multiply\count231 by \count221 
\multiply\count232 by \count222 
\multiply\count233 by \count223 
\advance\count153 by \count231
\advance\count153 by \count232
\advance\count153 by \count233
\divide\count153 by 5
\fi%
\count211=\count141
\count212=\count142
\count213=\count143
\count221=\count151
\count222=\count152
\count223=\count153
}

\def\serr#1#2{\power{3}{#1}
\loop 
\advance\count91 by -1 
\count92=#1 \count95=\count91 
\count211=0
\count212=62500
\count213=125000
\count221=0
\count222=109375
\count223=0
{\loop 
\advance\count92 by -1 
\count93=\count95 
\divide\count95 by 3 
\count94=\count95 
\multiply\count94 by 3 
\advance\count93 by -\count94 
\trooh{\count93} 
\ifnum \count92>0 \repeat \trigg{\count211}{\count221}{\count212}{\count222}{\count213}{\count223}{#2}}
\ifnum \count91>0 \repeat} 

\BBF{figSig}
{\Sig.}
{\begin{picture}(96,100)(-0,0)\setlength{\unitlength}{0.5pt}\small
\setlength{\unitlength}{3.2pt}
\put(-2.5,-2.5){{$v_1$}} 
\put(13.75,28.75){{$v_2$}} 
\put(30.625,-2.5){{$v_3$}} 
\setlength{\unitlength}{.006pt} 
\SIGA{7}{\Tiny.}
\end{picture}}

One can see that the analysis on the \Sig\ 
is analogous in many respects to the analysis of the fractal \Lp\ on an 
interval we described in \Sec{s4}. This analysis 
does not depend on 
the particular embedding of 
the \Sig\ in $\mathbb R^2$, and so a particular 
choice of the points $v_j$ does not matter as long as they are not 
collinear. 

\BDF{df19}{
 For each $n \geqslant 0$ we define $V_\ii n$ by 
$$
V_\ii {n\ps+1} =
F_\ii 1(V_\ii {n})\bigcup F_\ii 2(V_\ii {n})
\bigcup F_\ii 3(V_\ii {n})
$$
where 
$V_\ii 0=\{v_1,v_2,v_3\}=\partial S$ 
is the set of vertices of $S$. 

}


\def\Tri#1#2#3#4{{
\put(#1,#2){\line(3,5){#3}}
\put(#1,#2){\line(3,0){#4}}
\count223=#1
\advance\count223 by #4
\put(\count223,#2){\line(-3,5){#3}}}} 

\def\Trii#1#2#3#4#5{{
\put(#1,#2){\circle*{#5}}
\count223=#1
\advance\count223 by #4
\put(\count223,#2){\circle*{#5}}
\count226=#1
\advance\count226 by #3
\count224=#3
\divide\count224 by 3
\multiply\count224 by 5
\count225=#2
\advance\count225 by \count224
\put(\count226,\count225){\circle*{#5}}
}}

\def\Triii#1#2#3#4#5{{
\put(#1,#2){\circle{#5}}
\count223=#1
\advance\count223 by #4
\put(\count223,#2){\circle{#5}}
\count226=#1
\advance\count226 by #3
\count224=#3
\divide\count224 by 3
\multiply\count224 by 5
\count225=#2
\advance\count225 by \count224
\put(\count226,\count225){\circle{#5}}
}}

\def\Spic#1#2#3#4{{
\count205=#2\count206=#2\count207=#2\count208=#2\count209=#2
\count210=#1\count211=#3\count214=#3\count212=#4
\divide\count210 by 2
{\ifnum\count210>0
\Spic{\count210}{#2}{#3}{#4}
\multiply\count207 by 3
\multiply\count208 by 6
\multiply\count209 by 5
\multiply\count207 by \count210
\multiply\count208 by \count210
\multiply\count209 by \count210
{\advance\count211 by \count207
\advance\count212 by \count209
\Spic{\count210}{#2}{\count211}{\count212}}
{\advance\count214 by \count208
\Spic{\count210}{#2}{\count214}{#4}}
\else
\multiply\count205 by 3
\multiply\count206 by 6
\Tri{#3}{#4}{\count205}{\count206}
\fi
}}}

\def\Spicc#1#2#3#4#5{{
\count205=#2\count206=#2\count207=#2\count208=#2\count209=#2
\count210=#1\count211=#3\count214=#3\count212=#4
\divide\count210 by 2
{\ifnum\count210>0
\Spicc{\count210}{#2}{#3}{#4}{#5}
\multiply\count207 by 3
\multiply\count208 by 6
\multiply\count209 by 5
\multiply\count207 by \count210
\multiply\count208 by \count210
\multiply\count209 by \count210
{\advance\count211 by \count207
\advance\count212 by \count209
\Spicc{\count210}{#2}{\count211}{\count212}{#5}}
{\advance\count214 by \count208
\Spicc{\count210}{#2}{\count214}{#4}{#5}}
\else
\multiply\count205 by 3
\multiply\count206 by 6
\Trii{#3}{#4}{\count205}{\count206}{#5}
\fi
}}}

\def\Spiccc#1#2#3#4#5{{
\count205=#2\count206=#2\count207=#2\count208=#2\count209=#2
\count210=#1\count211=#3\count214=#3\count212=#4
\divide\count210 by 2
{\ifnum\count210>0
\Spiccc{\count210}{#2}{#3}{#4}{#5}
\multiply\count207 by 3
\multiply\count208 by 6
\multiply\count209 by 5
\multiply\count207 by \count210
\multiply\count208 by \count210
\multiply\count209 by \count210
{\advance\count211 by \count207
\advance\count212 by \count209
\Spiccc{\count210}{#2}{\count211}{\count212}{#5}}
{\advance\count214 by \count208
\Spiccc{\count210}{#2}{\count214}{#4}{#5}}
\else
\multiply\count205 by 3
\multiply\count206 by 6
\Triii{#3}{#4}{\count205}{\count206}{#5}
\fi
}}}

\BBF{figDASig}
{Discrete approximations to the \Sig.}
{%
\begin{picture}(60,50)(0,0) \setlength{\unitlength}{0.5pt}\small\thinlines
\put(-20,70){$V_\ii {0}:$}
{\def\LLL{6pt}
\setlength{\unitlength}{10pt}
\Spicc{1}{1}{0}{0}{.5}%
\Spic{1}{1}{0}{0}%
}
\end{picture}\ \ \ 
\ \ \ 
\ \ \ 
\begin{picture}(60,50)(0,0) \setlength{\unitlength}{0.5pt}\small\thinlines
\put(-20,70){$V_\ii {1}:$}
{\def\LLL{6pt}
\setlength{\unitlength}{5pt}
\Spicc{2}{1}{0}{0}{.9}%
\Spic{2}{1}{0}{0}%
}
\end{picture}\ \ \ 
\ \ \ 
\ \ \ 
\begin{picture}(60,50)(0,0) \setlength{\unitlength}{0.5pt}\small\thinlines
\put(-20,70){$V_\ii {2}:$}
{\def\LLL{6pt}
\setlength{\unitlength}{2.5pt}
\Spicc{4}{1}{0}{0}{1.5}%
\Spic{4}{1}{0}{0}%
}
\end{picture}%
}

On each $V_n$ there is a notion of neighboring points, 
giving $V_n$ a natural graph structure, 
which can be defined inductively. 
In \Fig{figDASig} the neighboring points are connected 
by line segments, and in what follows are denoted by the symbol 
``$x\sim y$''. 


\BDF{df20}{ The \bydef{discrete Dirichlet (energy) form} on $V_n$ is 
 $$
\mathcal E_n (\ddD{f}) = 
\big(\tfrac 53\big)^ n 
{\sum\limits_ {\genfrac{}{}{0pt}{2}{x,y\in V_\ii {n}}{y \sim x}}}
(f\sm{(y)} {-} f\x )^\ii 2
$$ 
and the \bydef{Dirichlet (energy) form} on $I$ is 
$$
\mathcal E (\ddD{f}) = \lim\limits_\ii {n{\to}\infty} 
\mathcal E_n (\ddD{f}) 
$$ 
if this limit exist and is finite. 
By the next proposition this limit is increasing and so it always exists as a nonnegative real number or $+\infty$. 

A function $h$ is \bydef{harmonic} if it minimizes the 
energy given the three boundary values.} 

\BPROP{prop22}{We have 
$
\mathcal E_{n{+}1} (\ddD{f})\geqslant\mathcal E_n (\ddD{f})
$ 
 for any function $f$, and 
$$
\mathcal E_{n{+}1}(\ddD{h})= 
\mathcal E_n (\ddD{h})=\mathcal E(\ddD{h})
$$
for a harmonic $h$.} 

\BDF{df23}{The \bydef{discrete \Lp s} on $V_n$ are 
$$ 
\Delta_n f\x = \tfrac14 
{\sum\limits_ {\genfrac{}{}{0pt}{2}{y\in V_\ii {n}}{y \sim x}}}
 f\sm{(y)} {-} f\x , \quad x{\in}V_n\backslash V_0
$$ 
and the \Lp\ on $S$ is 
$$
\Delta_\mu f\x = \lim\limits_\ii {n{\to}\infty} 5^ n 
\Delta_n f\x 
$$ 
if this limit exists and $\Delta_\mu f$ is continuous.}

The next theorem is essentially due to J. Kigami \cite{Ki0,Ki1,Ki}, 
although many relevant results are contained in 
\cite{BP,Fu,FS,Go,KL,Ko}. 
\BTHMm{}{
$\mathcal E$ is a local regular Dirichlet form on $S$ which is 
self-similar in the sense that
$$ 
\mathcal E(\ddD{f})=
\tfrac53\sum_{\sm{\sm{j=1,2,3}}}\mathcal E(\ddD{f\sm{\circ}F_j}).
$$
In particular, the domains of $\mathcal E$ and $\Delta_\mu $ 
are contained and dense in the space of continuous functions. 
The \Lp\ satisfies the 
\bydef{Gauss--Green (integration by parts) formula}: 
$$
 \mathcal E (\ddD{f})=-6\int_S f\Delta_\mu f d\mu + 
\sum_{p\in\partial S} f(p)\partial f(p) 
$$
where $\mu$ is the normalized Hausdorff measure, 
which is self-similar with weights $\frac13,\frac13,\frac13$ 
$$ 
\mu =\tfrac13\sum_{\sm{\sm{j=1,2,3}}} \mu \circ F_j,
$$
and $\partial f(p)$ is a certain derivative at the boundary. 
The (Dirichlet, Neumann) \Lp\ has discrete spectrum and compact 
resolvent. Green's function 
$g(x,y)$ is continuous. 
If $\rho(\lambda)$ is the eigenvalue counting function 
of the \Lp\ $\Delta_\mu$, then
$$
0<
\liminf_{\lambda{\to}\infty}
\frac{\rho(\lambda)}{\lambda^{d_s/2}}
< 
\limsup_{\lambda{\to}\infty}
\frac{\rho(\lambda)}{\lambda^{d_s/2}}
<\infty
$$
where the spectral dimension is 
$$
1<d_s=\tfrac{\log9}{\log5}<2.
$$ 
} 

Next, we present two more theorems illustrating the ``fractal'' 
properties of the \Lp\ $\Delta_\mu$. 

\BTHMm{\cite{BP,MS}}{The $\mu$--heat kernel of the \Lp\ $\Delta_\mu$
(that is, the transition density $p^\mu_t(x,y)$ of the $\mu$--diffusion)
has the non-Gaussian ``fractal'' asymptotics
$$
 p^\mu_t(x,y) \asymp
C\,\, t^{-\beta} 
\exp\left\{-C\,{t}^{-\alpha}\,{d_\ii{S}(x,y)^\gamma}\right\}
$$
where\, 
$\alpha{=}\frac{\log2}{\log5-\log2}$,\,
$\beta {=}\frac{\log3}{\log5}$\, and\, 
$\gamma{=}\frac{\log5}{\log5-\log2}$.}

\BTHMm{\cite{BST}}{
If $f\in \mathcal Dom\Delta_\mu$ then 
$f^ 2\notin \mathcal Dom\Delta_\mu$ unless $f$ is constant.}

Our main result in this section is the following theorem. 

\BTHM{thm25}{The zeta function of the \Lp\ on 
the \Sig\ is given by 
\begin{multline}\LL{eeeee}
\zeta{\ind{{\Delta_\mu}}}(s)=
\zetaa{\frac34}   {R}(s)
\mbox{$\cdot\frac{5^{- {s} /2}}{2}\pms\cdot
\left(
\frac{1}{1-3\cdot5^{- {s} /2}}+\frac{3}{1-5^{- {s} /2}}
\right)$}
+\\
\zetaa{\frac54}   {R}(s)
\mbox{$\cdot\frac{5^{- {s} }}{2}\pms\cdot
\left(
\frac{3}{1-3\cdot5^{- {s} /2}}-\frac{1}{1-5^{- {s} /2}}
\right)$}\end{multline} 
where $R(z)=z(5-4z)$. 
In particular, 
$\zeta{\ind{{\Delta_\mu}}}(s)$ has a meromorphic continuation for 
$ \text{Re}(s)\pms> -\varepsilon$ where $\varepsilon\geqslant1$. 
The poles are contained in 
$${\left\{\frac{2in\pi}{\log 5},\frac{\log 9+2in\pi}{\log 5}: 
n\in\mathbb Z\right\}}.$$ 
We have 
$0<\frac{\log4}{\log5}=d_R<1<d_s=\frac{\log9}{\log5}<2$.}

\Pr 
By \cite{FS,T}, every eigenvalue of $\Delta_\mu$ has the form 
$$
\lambda=-5^m\lim\limits_{n{\to}\infty}5^n R^{\sm{\sm{\sm{-}}}n}(z_0)
$$ 
where 
$R^{\sm{\sm{\sm{-}}}n}(z_0)$ is a preimage of 
$z_0=3/4,5/4$ under the $n$-th 
iteration power of the polynomial $R(z)$. 
The multiplicity of such an eigenvalue is of the form 
$C_13^m+C_2$ where $C_1$ and $C_2$ are computed explicitly in \cite{FS,T}. 
More precisely, by \cite[Proposition 3.12]{T}, 
for the discrete \Lp\ we have that the multiplicity of $z\in R^{-m}\left(\frac34\right) $ 
is $\frac{3^{n-m-1}+3}2$ for $n\geqslant1$ and $0\leqslant m \leqslant n-1$, 
and the multiplicity of $z\in R^{-m}\left(\frac54\right) $ is $\frac{3^{n-m-1}-1}2$ for 
$n\geqslant 2$ and $0\leqslant m\leqslant n-2 $. Then the proof follows directly from \Thm{thmPoly}  and the definition of the \zfoap. 
 \rP

\BBF{figSigCSD}
{Complex spectral dimensions of the \Lp\ on the \Sig.}
{\begin{picture}(250,200)(-125,-100) 
\setlength{\unitlength}{0.5pt}\small
\thicklines\setlength{\unitlength}{0.4pt}
\put(-250,0){\vector(1,0){500}}
\put(0,-200){\vector(0,1){400}}
\put(5,5){{$0$}}
\put(103,7){$\sm{d_R}$}
\put(152,7){${1}$}
\put(-145,7){${-1}$}
\put(150,-5){\line(0,1){10}}
\put(-150,-5){\line(0,1){10}}
\put(213,7){${d_s}$}
\multiput(100,-180)(0,60){7}{\circle*{8}}
\multiput(210,-180)(0,60){7}{\circle{7}}
\multiput(0,-180)(0,60){7}{\circle{7}}
\multiput(-150,-200)(0,20){20}{\line(0,1){10}}
\thinlines
\multiput(-150,-150)(0,4){87}{\line(-5,-3){120}}
\end{picture}}

\BREM{rem26}{We have that the zeros of the ``geometric part'' cancel 
all the poles of the \zfotp in the product formula 
(\ref{eeeee}). 
We conjecture this also for other symmetric p.c.f. fractals \cite{MT}, except 
one dimensional. 
In \Fig{figSigCSD} the series of the canceled poles is denoted by black dots, 
and the other two series of poles are denoted by small circles.}

\BREM{rem26-}{Note that the \zfotp in (\ref{eeeee}) also describes 
the distribution of the eigenvalues corresponding to the 
non localized eigenfunctions (see \cite{BK,T}).}

\BREM{rem26--}{The method of Theorem \ref{thm25} applies to a large class of fractals 
(see, for instance, \cite{KrTe03,MT,T,T04,T01,T02}), and possibly to \cite{S}.}

\end{document}